%
% On the Frobenius integrability of certain holomorphic $p$-forms
% April 2000
%
% Jean-Pierre Demailly,
% Universit\'e de Grenoble I, Institut Fourier
% 38402 Saint-Martin d'H\`eres, France
%
%
% Plain-TeX file

% page setting

\magnification=1200
\pretolerance=500 \tolerance=1000 \brokenpenalty=5000
\hsize=12.5cm   
\vsize=19cm
\hoffset=0.4cm
\voffset=1cm   
\parskip3pt plus 1pt
\parindent=0.6cm
\let\sl=\it
\def\\{\hfil\break}

% new fonts definitions

\font\seventeenbf=cmbx10 at 17.28pt

\font\twelvebf=cmbx10 at 12pt
\font\eightbf=cmbx8
\font\sixbf=cmbx6

\font\eighti=cmmi8
\font\sixi=cmmi6

\font\eightrm=cmr8
\font\sixrm=cmr6

\font\eightsy=cmsy8
\font\sixsy=cmsy6

\font\eightit=cmti8
\font\eighttt=cmtt8
\font\eightsl=cmsl8

\font\seventeenbsy=cmbsy10 at 17.28pt

\font\twelvebsy=cmbsy10 at 12pt
\font\tenbsy=cmbsy10
\font\eightbsy=cmbsy8
\font\sevenbsy=cmbsy7
\font\sixbsy=cmbsy6
\font\fivebsy=cmbsy5

\font\tenmsa=msam10

\font\sevenmsa=msam7
\font\fivemsa=msam5
\newfam\msafam
  \textfont\msafam=\tenmsa
  \scriptfont\msafam=\sevenmsa
  \scriptscriptfont\msafam=\fivemsa

\font\tenmsb=msbm10
\font\eightmsb=msbm8
\font\sevenmsb=msbm7
\font\fivemsb=msbm5
\newfam\msbfam
  \textfont\msbfam=\tenmsb
  \scriptfont\msbfam=\sevenmsb
  \scriptscriptfont\msbfam=\fivemsb
\def\Bbb{\fam\msbfam\tenmsb}

\font\tenCal=eusm10
\font\sevenCal=eusm7
\font\fiveCal=eusm5
\newfam\Calfam
  \textfont\Calfam=\tenCal
  \scriptfont\Calfam=\sevenCal
  \scriptscriptfont\Calfam=\fiveCal
\def\Cal{\fam\Calfam\tenCal}

\font\teneuf=eusm10
\font\teneuf=eufm10
\font\seveneuf=eufm7
\font\fiveeuf=eufm5
\newfam\euffam
  \textfont\euffam=\teneuf
  \scriptfont\euffam=\seveneuf
  \scriptscriptfont\euffam=\fiveeuf

\font\seventeenbfit=cmmib10 at 17.28pt

\font\twelvebfit=cmmib10 at 12pt
\font\tenbfit=cmmib10
\font\eightbfit=cmmib8
\font\sevenbfit=cmmib7
\font\sixbfit=cmmib6
\font\fivebfit=cmmib5
\newfam\bfitfam
  \textfont\bfitfam=\tenbfit
  \scriptfont\bfitfam=\sevenbfit
  \scriptscriptfont\bfitfam=\fivebfit

% changing font sizes

\catcode`\@=11
\def\eightpoint{%
  \textfont0=\eightrm \scriptfont0=\sixrm \scriptscriptfont0=\fiverm
  \def\rm{\fam\z@\eightrm}%
  \textfont1=\eighti \scriptfont1=\sixi \scriptscriptfont1=\fivei
  \def\oldstyle{\fam\@ne\eighti}%
  \textfont2=\eightsy \scriptfont2=\sixsy \scriptscriptfont2=\fivesy
  \textfont\itfam=\eightit
  \def\it{\fam\itfam\eightit}%
  \textfont\slfam=\eightsl
  \def\sl{\fam\slfam\eightsl}%
  \textfont\bffam=\eightbf \scriptfont\bffam=\sixbf
  \scriptscriptfont\bffam=\fivebf
  \def\bf{\fam\bffam\eightbf}%
  \textfont\ttfam=\eighttt
  \def\tt{\fam\ttfam\eighttt}%
  \textfont\msbfam=\eightmsb
  \def\Bbb{\fam\msbfam\eightmsb}%
  \abovedisplayskip=9pt plus 2pt minus 6pt
  \abovedisplayshortskip=0pt plus 2pt
  \belowdisplayskip=9pt plus 2pt minus 6pt
  \belowdisplayshortskip=5pt plus 2pt minus 3pt
  \smallskipamount=2pt plus 1pt minus 1pt
  \medskipamount=4pt plus 2pt minus 1pt
  \bigskipamount=9pt plus 3pt minus 3pt
  \normalbaselineskip=9pt
  \setbox\strutbox=\hbox{\vrule height7pt depth2pt width0pt}%
  \let\bigf@ntpc=\eightrm \let\smallf@ntpc=\sixrm
  \normalbaselines\rm}
\catcode`\@=12

\def\eightpointbf{%
 \textfont0=\eightbf   \scriptfont0=\sixbf   \scriptscriptfont0=\fivebf
 \textfont1=\eightbfit \scriptfont1=\sixbfit \scriptscriptfont1=\fivebfit
 \textfont2=\eightbsy  \scriptfont2=\sixbsy  \scriptscriptfont2=\fivebsy
 \eightbf
 \baselineskip=10pt}

\def\tenpointbf{%
 \textfont0=\tenbf   \scriptfont0=\sevenbf   \scriptscriptfont0=\fivebf
 \textfont1=\tenbfit \scriptfont1=\sevenbfit \scriptscriptfont1=\fivebfit
 \textfont2=\tenbsy  \scriptfont2=\sevenbsy  \scriptscriptfont2=\fivebsy
 \tenbf}
        
\def\twelvepointbf{%
 \textfont0=\twelvebf   \scriptfont0=\eightbf   \scriptscriptfont0=\sixbf
 \textfont1=\twelvebfit \scriptfont1=\eightbfit \scriptscriptfont1=\sixbfit
 \textfont2=\twelvebsy  \scriptfont2=\eightbsy  \scriptscriptfont2=\sixbsy
 \twelvebf
 \baselineskip=14.4pt}

\def\seventeenpointbf{%
 \textfont0=\seventeenbf  \scriptfont0=\twelvebf  \scriptscriptfont0=\eightbf
 \textfont1=\seventeenbfit\scriptfont1=\twelvebfit\scriptscriptfont1=\eightbfit
 \textfont2=\seventeenbsy \scriptfont2=\twelvebsy \scriptscriptfont2=\eightbsy
 \seventeenbf
 \baselineskip=20.736pt}
 
% main item macros

\newdimen\srdim \srdim=\hsize
\newdimen\irdim \irdim=\hsize
\def\NOSECTREF#1{\noindent\hbox to \srdim{\null\dotfill ???(#1)}}
\def\SECTREF#1{\noindent\hbox to \srdim{\csname REF\romannumeral#1\endcsname}}
\def\INDREF#1{\noindent\hbox to \irdim{\csname IND\romannumeral#1\endcsname}}
\newlinechar=`\^^J
\def\openauxfile{
  \immediate\openin1\jobname.aux
  \ifeof1
  \message{^^JCAUTION\string: you MUST run TeX a second time^^J}
  \let\sectref=\NOSECTREF \let\indref=\NOSECTREF
  \else
  \input \jobname.aux
  \message{^^JCAUTION\string: if the file has just been modified you may 
    have to run TeX twice^^J}
  \let\sectref=\SECTREF \let\indref=\INDREF
  \fi
  \message{to get correct page numbers displayed in Contents or Index 
    Tables^^J}
  \immediate\openout1=\jobname.aux
  \let\END=\end \def\end{\immediate\closeout1\END}}
        
\newbox\titlebox   \setbox\titlebox\hbox{\hfil}
\newbox\sectionbox \setbox\sectionbox\hbox{\hfil}
\def\folio{\ifnum\pageno=1 \hfil \else \ifodd\pageno
           \hfil {\eightpoint\copy\sectionbox\kern8mm\number\pageno}\else
           {\eightpoint\number\pageno\kern8mm\copy\titlebox}\hfil \fi\fi}
\footline={\hfil}
\headline={\folio}

\def\titlerunning#1{\setbox\titlebox\hbox{\eightpoint #1}}
\def\title#1{\noindent\hfil$\smash{\hbox{\seventeenpointbf #1}}$\hfil
             \titlerunning{#1}\medskip}

\newcount\numbersection \numbersection=-1
\def\sectionrunning#1{\setbox\sectionbox\hbox{\eightpoint #1}
  \immediate\write1{\string\def \string\REF 
      \romannumeral\numbersection \string{%
      \noexpand#1 \string\dotfill \space \number\pageno \string}}}
\def\section#1{%
  \par\vskip0.666cm\penalty -100
  \vbox{\baselineskip=14.4pt\noindent{{\twelvepointbf #1}}}
  \vskip2pt
  \penalty 500
  \advance\numbersection by 1
  \sectionrunning{#1}}

\def\subsection#1|{%
  \par\vskip0.5cm\penalty -100
  \vbox{\noindent{{\tenpointbf #1}}}
  \vskip1pt
  \penalty 500}

\newcount\numberindex \numberindex=0  
\def\index#1#2{%
  \advance\numberindex by 1
  \immediate\write1{\string\def \string\IND #1%
     \romannumeral\numberindex \string{%
     \noexpand#2 \string\dotfill \space \string\S \number\numbersection, 
     p.\string\ \space\number\pageno \string}}}

\newdimen\itemindent \itemindent=\parindent

\def\item#1{\par\noindent\hangindent\itemindent%
            \rlap{#1}\kern\itemindent\ignorespaces}
\def\itemitem#1{\par\noindent\hangindent2\itemindent%
            \kern\itemindent\rlap{#1}\kern\itemindent\ignorespaces}
\def\itemitemitem#1{\par\noindent\hangindent3\itemindent%
            \kern2\itemindent\rlap{#1}\kern\itemindent\ignorespaces}

\long\def\claim#1|#2\endclaim{\par\vskip 5pt\noindent 
{\tenpointbf #1.}\ {\sl #2}\par\vskip 5pt}

\def\today{\ifcase\month\or
January\or February\or March\or April\or May\or June\or July\or August\or
September\or October\or November\or December\fi \space\number\day,
\number\year}

\catcode`\@=11
\newcount\@tempcnta \newcount\@tempcntb 
\def\timeofday{{%
\@tempcnta=\time \divide\@tempcnta by 60 \@tempcntb=\@tempcnta
\multiply\@tempcntb by -60 \advance\@tempcntb by \time
\ifnum\@tempcntb > 9 \number\@tempcnta:\number\@tempcntb
  \else\number\@tempcnta:0\number\@tempcntb\fi}}
\catcode`\@=12

\def\bibitem#1&#2&#3&#4&%
{\hangindent=1.66cm\hangafter=1
\noindent\rlap{\hbox{\eightpointbf #1}}\kern1.66cm{\rm #2}{\sl #3}{\rm #4.}} 

% blackboard symbols

\def\bC{{\Bbb C}}

\def\bR{{\Bbb R}}

% gothic symbols

% calligraphic symbols
\def\cO{{\Cal O}}
\def\cS{{\Cal S}}

% special symbols

\def\ld{,\ldots,}

\def\square{{\hfill \hbox{
\vrule height 1.453ex  width 0.093ex  depth 0ex
\vrule height 1.5ex  width 1.3ex  depth -1.407ex\kern-0.1ex
\vrule height 1.453ex  width 0.093ex  depth 0ex\kern-1.35ex
\vrule height 0.093ex  width 1.3ex  depth 0ex}}}

\def\hexnbr#1{\ifnum#1<10 \number#1\else
 \ifnum#1=10 A\else\ifnum#1=11 B\else\ifnum#1=12 C\else
 \ifnum#1=13 D\else\ifnum#1=14 E\else\ifnum#1=15 F\fi\fi\fi\fi\fi\fi\fi}
\def\msatype{\hexnbr\msafam}
\def\msbtype{\hexnbr\msbfam}
\mathchardef\restriction="3\msatype16   
\mathchardef\compact="3\msatype62
\mathchardef\smallsetminus="2\msbtype72   \let\ssm\smallsetminus
\mathchardef\subsetneq="3\msbtype28
\mathchardef\supsetneq="3\msbtype29
\mathchardef\leqslant="3\msatype36   \let\le\leqslant
\mathchardef\geqslant="3\msatype3E   \let\ge\geqslant
\mathchardef\ltimes="2\msbtype6E
\mathchardef\rtimes="2\msbtype6F

% hats and tildes and over/underlines

\let\wh=\widehat
\let\text=\hbox
\def\build#1|#2|#3|{\mathrel{\mathop{\null#1}\limits^{#2}_{#3}}}

% mathematical operators

\def\del{{\partial}}
\def\dbar{{\overline\partial}}
\def\ddbar{{\partial\overline\partial}}

% subscript operands

% figures inserted as PostScript files
%\special{header=/home/demailly/grlib.ps}
\long\def\InsertFig#1 #2 #3 #4\EndFig{\par
\hbox{\hskip #1mm$\vbox to#2mm{\vfil\special{" 
(/home/demailly/psinputs/grlib.ps) run
#3}}#4$}}
\long\def\LabelTeX#1 #2 #3\ELTX{\rlap{\kern#1mm\raise#2mm\hbox{#3}}}

% main text

\openauxfile

\title{On the Frobenius integrability}
\title{of certain holomorphic $p$-forms}
\titlerunning{Frobenius integrability of certain holomorphic $p$-forms}
\medskip
\centerline{\twelvebf Jean-Pierre Demailly}
\vskip25pt
\line{\hfill \it Dedicated to Professor Hans Grauert, on the occasion
of his 70th birthday} 
\vskip20pt

\noindent
{\bf Abstract.}
The goal of this note is to exhibit the integrability properties (in
the sense of the Frobenius theorem) of holomorphic $p$-forms with
values in certain line bundles with semi-negative curvature on a compact
K\"ahler manifold. There are in fact very strong restrictions, both on
the holomorphic form and on the curvature of the semi-negative line bundle. 
In particular, these observations provide interesting information on the
structure of projective manifolds which admit a contact structure:
either they are Fano manifolds or, thanks to results of
Kebekus-Peternell-Sommese-Wisniewski, they are biholomorphic to the
projectivization of the cotangent bundle of another suitable
projective manifold.

\section{1. Main results}

Recall that a holomorphic line bundle $L$ on a compact complex
manifold is said to be {\it pseudo-effective} if $c_1(L)$ contains a closed
positive $(1,1)$-current $T$, or equivalently, if $L$ possesses a
(possibly singular) hermitian metric $h$ such that the curvature current
$T=\Theta_h(L)=-i\ddbar\log h$ is nonnegative. If $X$ is projective,
$L$ is pseudo-effective if and only if $c_1(L)$ belongs to the closed
cone of $H^{1,1}_{\bR}(X)$ generated by classes of effective divisors
(see [Dem90, 92]). Our main result is

\claim Main Theorem|Let $X$ be a compact K\"ahler manifold. Assume that
there exists a pseudo-effective line bundle $L$ on $X$ and a nonzero
holomorphic section $\theta\in H^0(X,\Omega^p_X\otimes L^{-1})$, where
$0\le p\le n=\dim X$. Let $\cS_\theta$ be the coherent subsheaf of germs 
of vector fields $\xi$ in the tangent sheaf $T_X$, such that the contraction
$i_\xi\theta$ vanishes. Then $\cS_\theta$ is integrable, namely 
$[\cS_\theta,\cS_\theta]\subset\cS_\theta$, and $L$ has flat curvature
along the leaves of the $($possibly singular$)$ foliation defined 
by~$\cS_\theta$.
\endclaim

Before entering into the proof, we discuss several consequences. If
$p=0$ or $p=n$, the result is trivial (with $\cS_\theta=T_X$ and
$\cS_\theta=0$, respectively). The most interesting case is $p=1$.

\claim Corollary 1|In the above situation, if the line bundle 
$L\to X$ is pseudo-effective and
\hbox{$\theta\in H^0(X,\Omega^1_X\otimes L^{-1})$} is a nonzero section, 
the subsheaf $\cS_\theta$ defines a holomorphic foliation of codimension 
$1$ in $X$, that is, $\theta\wedge d\theta=0$.
\endclaim

We now concentrate ourselves on the case when $X$ is a {\it contact 
manifold}, i.e.\ $\dim X=n=2m+1$, $m\ge 1$, and there exists a form
$\theta\in H^0(X,\Omega^1_X\otimes L^{-1})$, called the {\it contact form},
such that $\theta\wedge(d\theta)^m\in H^0(X,K_X\otimes L^{-m-1})$
has no zeroes. Then $\cS_\theta$ is a codimension $1$ locally free
subsheaf of $T_X$ and there are dual exact sequences
$$
0\to L\to\Omega^1_X\to\cS_\theta^\star\to 0,\qquad
0\to \cS_\theta\to T_X\to L^\star\to 0.
$$
The subsheaf $\cS_\theta\subset T_X$ is said to be the 
{\it contact structure} of $X$. The assumption that $\theta\wedge(d\theta)^m$
does not vanish implies that $K_X\simeq L^{m+1}$. In that case, the subsheaf
is not integrable, hence $L$ and $K_X$ cannot be pseudo-effective.

\claim Corollary 2|If $X$ is a compact K\"ahler manifold admitting a
contact structure, then $K_X$ is not pseudo-effective, in particular
the Kodaira dimension $\kappa(X)$ is equal to $-\infty$.
\endclaim

The fact that $\kappa(X)=-\infty$ had been observed previously by
St\'ephane Druel [Dru98]. In the projective context, the minimal model
conjecture would imply (among many other things) that the conditions
$\kappa(X)=-\infty$ and ``$K_X$ non pseudo-effective'' are equivalent,
but a priori the latter property is much stronger (and in large
dimensions, the minimal model conjecture still seems far beyond
reach!)

\claim Corollary 3|If $X$ is a compact K\"ahler manifold with a
contact structure and with second Betti number $b_2=1$, then
$K_X$ is negative, i.e., $X$ is a Fano manifold.
\endclaim

Actually the Kodaira embedding theorem shows that the K\"ahler
manifold $X$ is projective if $b_2=1$, and then every line bundle is
either positive, flat or negative. As $K_X$ is not pseudo-effective it
must therefore be negative. In that direction, Boothby [Boo61], Wolf
[Wol65] and Beauville [Bea98] have exhibited a natural construction of
contact Fano manifolds. Each of the known examples is obtained as a
homogeneous variety which is the unique closed orbit in the
projectivized (co)adjoint representation of a simple algebraic Lie
group. Beauville's work ([Bea98], [Bea99]) provides strong evidence
that this is the complete classification in the case $b_2=1$.

We now come to the case $b_2\ge 2$. If $Y$ is an arbitrary compact 
K\"ahler manifold, the bundle
$X=P(T_Y^\star)$ of hyperplanes of $T_Y$ has a contact structure
associated with the line bundle $L=\cO_X(-1)$. Actually, if 
$\pi:X\to Y$ is the canonical projection, one can define a contact 
form $\theta\in H^0(X,\Omega^1_X\otimes L^{-1})$ by setting
$$
\theta(x)=\theta(y,[\xi])=\xi^{-1}\pi^\star\xi=\xi^{-1}\sum_{1\le j\le p}
\xi_jdy_j,\qquad p=\dim Y,
$$
at every point $x=(y,[\xi])\in X$, $\xi\in T_{Y,y}^\star\ssm\{0\}$ 
(observe that $\xi\in L_x=\cO_X(-1)_x$). Morever $b_2(X)=1+b_2(Y)\ge 2$. 
Conversely, Kebekus, Peternell, Sommese and Wi{\'s}niewski [KPSW] have 
recently shown that every projective algebraic manifold~$X$ such that
\item{(i)} $X$ has a contact structure,
\item{(ii)} $b_2\ge 2$,
\item{(iii)}  $K_X$ is not nef (numerically effective)
\smallskip\noindent
is of the form $X=P(T_Y^\star)$ for some projective algebraic manifold $Y$.
However, the condition that $K_X$ is not nef is implied by the fact that
$K_X$ is not pseudo-effective. Hence we get

\claim Corollary 4|If $X$ is a contact projective manifold with
$b_2\ge 2$, then $X$ is a projectivized hyperplane bundle 
$X=P(T_Y^\star)$ associated with some projective manifold $Y$.
\endclaim

The K\"ahler case of corollary 4 is still unsolved, as the proof of [KPSW]
heavily relies on Mori theory (and, unfortunately, the extension of Mori
theory to compact K\"ahler manifolds remains to be settled~$\ldots$). 

I would like to thank Arnaud Beauville, Fr\'ed\'eric Campana, Stefan
Kebekus and Thomas Peternell for illuminating discussions on these
subjects. The present work was written during a visit at G\"ottingen
University, on the occasion of a colloquium in honor of Professor Hans
Grauert for his 70th birthday.

\section{2. Proof of the Main Theorem}

In some sense, the proof is just a straightforward integration by parts, 
but there are slight technical difficulties due to the fact that we have
to work with singular metrics.

Let $X$ be a compact K\"ahler manifold, $\omega$ the K\"ahler metric, and
let $L$ be a pseudo-effective line bundle on $X$. We select a
hermitian metric $h$ on $L$ with nonnegative curvature current
$\Theta_h(L)\ge 0$, and let $\varphi$ be the plurisubharmonic weight
of the metric $h$ in any local trivialisation $L_{|U}\simeq U\times\bC$.
In other words, we have
$$
\Vert\xi\Vert_h^2=|\xi|^2e^{-\varphi(x)},\qquad
\Vert\xi^\star\Vert_{h^\star}^2=|\xi^\star|^2e^{\varphi(x)}
$$
for all $x\in U$ and $\xi\in L_x$, $\xi^\star\in L^{-1}$. We then
have a Chern connection $\nabla=\del_{h^\star}+\dbar$ acting on all 
$(p,q)$-forms $f$ with values in $L^{-1}$, given locally by
$$
\del_\varphi f =e^{-\varphi}\del(e^\varphi f)=\del f+\del\varphi\wedge f
$$
in every trivialization $L_{|U}$. Now, assume that there is a holomorphic
section $\theta\in H^0(X,\Omega^p_X\otimes L^{-1})$, i.e., a 
$\dbar$-closed $(p,0)$ form $\theta$ with values in $L^{-1}$.
We compute the global $L^2$ norm
$$
\int_X \{\del_{h^\star}\theta,\del_{h^\star}\theta\}_{h^\star}
\wedge\omega^{n-p-1}=
\int_X e^\varphi\del_\varphi\theta\wedge
\overline{\del_\varphi\theta}\wedge\omega^{n-p-1}
$$
where $\{\kern5pt,\kern4pt\}_{h^\star}$ is the natural sesquilinear
pairing sending pairs of $L^{-1}$-valued forms of type $(p,q)$,
$(r,s)$ into $(p+s,q+r)$ complex valued forms. The right hand side is of
course only locally defined, but it explains better how the forms are
calculated, and also all local representatives glue together into a
well defined global form; we~will therefore use the latter notation as 
if it were global. As
$$
d\big(e^\varphi\theta\wedge
\overline{\del_\varphi\theta}\wedge\omega^{n-p-1}\big)=
e^\varphi\del_\varphi\theta\wedge
\overline{\del_\varphi\theta}\wedge\omega^{n-p-1}+(-1)^p
e^\varphi\theta\wedge\overline{\dbar\del_\varphi\theta}\wedge\omega^{n-p-1}
$$
and $\dbar\del_\varphi\theta=\dbar\del\varphi\wedge\theta$, an integration 
by parts via Stokes theorem yields
$$
\int_X e^\varphi\del_\varphi\theta\wedge
\overline{\del_\varphi\theta}\wedge\omega^{n-p-1}
=-(-1)^p\int_X e^\varphi\ddbar\varphi\wedge\theta\wedge\overline\theta
\wedge\omega^{n-p-1}.
$$
These calculations need a word of explanation, since $\varphi$ is in general
singular. However, it is well known that the $i\ddbar$ of a plurisubharmonic
function is a closed positive current, in particular
$$
i\ddbar(e^\varphi)=e^\varphi(i\del\varphi\wedge\dbar\varphi+i\ddbar\varphi)
$$
is positive and has measure coefficients. This shows that $\del\varphi$ is
$L^2$ with respect to the weight $e^\varphi$, and similarly that 
$e^\varphi\ddbar\varphi$ has locally finite measure coefficients.
Moreover, the results of [Dem92] imply that there is a decreasing sequence
of metrics $h^\star_\nu$ and corresponding weights $\varphi_\nu
\downarrow\varphi$, such that $\Theta_{h_\nu}\ge -C\omega$ with a fixed
constant $C>0$ (this claim is in fact much weaker than the results of
[Dem92], and easy to prove e.g.\ by using convolutions in suitable
coordinate patches and a standard gluing technique). Now, the results
of Bedford-Taylor [BT76, BT82] applied to the uniformly bounded functions
$e^{c\varphi_\nu}$, $c>0$, imply that we have local weak convergence
$$
e^{\varphi_\nu}\ddbar\varphi_\nu\to e^{\varphi}\ddbar\varphi,\quad
e^{\varphi_\nu}\del\varphi_\nu\to e^{\varphi}\del\varphi,\quad
e^{\varphi_\nu}\del\varphi_\nu\wedge\dbar\varphi_\nu\to 
e^{\varphi}\del\varphi\wedge\dbar\varphi,
$$
possibly afting adding $C'|z|^2$ to the $\varphi_\nu$'s to make them
plurisubharmonic. This is enough to justify the calculations. Now, we
take care of signs, using the fact that $i^{p^2}\theta\wedge\overline
\theta\ge 0$ whenever $\theta$ is a $(p,0)$-form. Our previous equality 
can be rewritten
$$ 
\int_X e^\varphi\,i^{(p+1)^2}\del_\varphi\theta\wedge
\overline{\del_\varphi\theta}\wedge\omega^{n-p-1}
=-\int_X e^\varphi\,i\ddbar\varphi\wedge i^{p^2}\theta\wedge\overline\theta
\wedge\omega^{n-p-1}.
$$
Since the left hand side is nonnegative and the right hand side is nonpositive,
we conclude that $\del_\varphi\theta=0$ almost everywhere, i.e.\
$\del\theta=-\del\varphi\wedge\theta$ almost everywhere. The formula
for the exterior derivative of a $p$-form reads
$$
\leqalignno{
d\theta(\xi_0\ld\xi_p)
&=\sum_{0\le j\le p}(-1)^j\xi_j\cdot \theta(\xi_0\ld\wh{\xi_j}\ld\xi_p)\cr
&\qquad{}
+\sum_{0\le j<k\le p}(-1)^{j+k}\theta([\xi_j,\xi_k],\xi_0\ld\wh{\xi_j},
\ldots,\wh{\xi_k}\ld\xi_p).&(\star)\cr}
$$
If two of the vector fields -- say $\xi_0$ and $\xi_1$ -- lie in
$\cS_\theta$, then
$$
d\theta(\xi_0\ld\xi_p)=-(\del\varphi\wedge\theta)(\xi_0\ld\xi_p)=0
$$
and all terms in the right hand side of $(\star)$ are also zero,
except perhaps the term $\theta([\xi_0,\xi_1],\xi_2\ld\xi_p)$. We
infer that this term must vanish. Since this is true for
arbitrary vector fields $\xi_2\ld\xi_p$, we conclude that
$[\xi_0,\xi_1]\in \cS_\theta$ and that $\cS_\theta$ is integrable.

The above arguments also yield strong restrictions on the
hermitian metric $h$. In fact the equality $\del\theta=-\del\varphi
\wedge\theta$ implies $\ddbar\varphi\wedge\theta=0$ by taking
the~$\dbar$.  Fix a smooth point in a leaf of the foliation, and local
coordinates $(z_1\ld z_n)$ such that the leaves are given by
$z_1=c_1\ld z_r=c_r$ ($c_i={}$constant), in a neighborhood of that
point. Then $\cS_\theta$ is generated by $\del/\del z_{r+1}\ld
\del/\del z_n$, and $\theta$ depends only on $dz_1\ld dz_r$. This
implies that $\del^2\varphi/\del z_j\del\overline z_k=0$ for
$j,k>r$, in other words $(L,h)$ has flat curvature along the leaves 
of the foliation. The main theorem is proved.

\bigskip\null
\section{References}
\medskip

{\eightpoint
  
\bibitem[Bea98]&Beauville, A.:& Fano contact manifolds and nilpotent orbits;&
  Comm.\ Math.\ Helv., {\bf 73}$\,$(4) (1998), 566--583&

\bibitem[Bea99]&Beauville, A.:& Riemannian holonomy and Algebraic Geometry;&
  Duke/alg-geom preprint 9902110, (1999)&

\bibitem[Boo61]&Boothby, W.:& Homogeneous complex contact manifolds;&
  Proc.\ Symp.\ Pure Math.\ {\bf 3} (Differential Geometry) (1961) 144--154&

\bibitem[BT76]&Bedford, E., Taylor, B.A.:& The Dirichlet problem for a complex
Monge-Amp\`ere equation;& Invent.\ Math.\ {\bf 37} (1976) 1--44&

\bibitem[BT82]&Bedford, E., Taylor, B.A.:& A new capacity for plurisubharmonic
functions;& Acta Math.\ {\bf 149} (1982) 1--41&

\bibitem[Dem90]&Demailly, J.-P.:& Singular hermitian metrics on
  positive line bundles;& Proceedings of the Bayreuth conference
  ``Complex algebraic varieties'', April~2-6, 1990, edited by
  K.~Hulek, T.~Peternell, M.~Schneider, F.~Schreyer, Lecture Notes in
  Math.\ n${}^\circ\,$1507, Springer-Verlag, 1992&

\bibitem[Dem92]&Demailly, J.-P.:& Regularization of closed positive currents 
  and Intersection Theory;& J.\ Alg.\ Geom.\ {\bf 1} (1992), 361--409&

\bibitem[Dru98]&Druel, S.:& Contact structures on algebraic $5$-dimensional 
  manifolds;& C.R.\ Acad.\ Sci.\ Paris, {\bf 327} (1998), 365--368&

\bibitem[KPSW]&Kebekus, S., Peternell, Th., Sommese, A.J.\ and 
  Wi{\'s}niewski, J.A.:& Projective Contact Manifolds;& preprint 1999,
  to appear in Inventiones Math&

\bibitem[Wol65]&Wolf, J.:& Complex homogeneous contact manifolds and
  quaternionic symmetric spaces;& J.\ Math.\ Mech.\ {\bf 14} (1965) 
  1033--1047&

}
\vskip20pt
\noindent
(version of April 13, 2000, printed on \today)
\vskip20pt

\parindent=0cm
Jean-Pierre Demailly\\
Universit\'e de Grenoble I,
D\'epartement de Math\'ematiques, Institut Fourier\\
38402 Saint-Martin d'H\`eres, France\\
{\it e-mail:\/} {\tt demailly@ujf-grenoble.fr}
\end